\documentclass[twocolumn]{article}
\title{A smooth variant of the Afriat- Varian
theorem}
\author{Alexander A. Shananin
\thanks{Supported by RHSF grant 01-02-00481a}
\\ Moscow State University
\\ Vorob'evy Gory
\\ Moscow, Russia.
\\ {\tt e-mail: shan@ccas.ru}
\and
Sergey P. Tarasov
\thanks{Supported in part by RHSF grant 01-02-00481a}
\\ Computing center of RAS \\ Vavilova 40 \\ 117967
Moscow GSP-1, Russia.\\
{\tt e-mail: sergey@ccas.ru}
}
\date{ }
  \newtheorem{theorem}{Theorem}
  \newtheorem{definition}{Definition}
  \newtheorem{lemma}{Lemma}
  \newtheorem{rK}{Remark}
  
  \newtheorem {pF}{Proof}

\def\RR{{\bf R}}

\def\F{{\cal F}}

\def\Q{{\cal Q}}
\def\P{{\cal P}}
\def\X{{\cal X}}
\def\BB{{\bf B}}

\begin{document}
\maketitle

\begin{abstract}
We present a simple geometric construction for smoothing polyhedral
utility functions.

\end{abstract}

\section{Introduction}

The problem of the integrability of the demand functions is studied in
mathematical economics for more than hundred years. (The functional
dependencies between the volume of purchases and the corresponding
prices are called the demand functions.) It seems that italian
economist J.Antonelli was one of the first to notice in 1886
(cf \cite{Pref})
that
under certain conditions the demand functions are rationalized, i.e.
they can be obtained via maximization of the utility function under
budget constraint. The attention to this problem was revived in 1905 in
connection with the discourse between V.Volterra and V.Pareto about the
interpretation of the demand functions' rationalization conditions
\cite{Pref}. V.Volterra pointed that Pareto's result on the
rationalizability of the demand functions for two kinds of
goods could not be generalized to the case of more than two kinds of
goods. However, it is known as well that rationalization is an
implicit assumption for the computation of the price indices and the
volume of purchases in economic statistics. V.Pareto has tried to
justify rationalization conditions and to show that they should be
always valid, but in 1915 E.Slutsky \cite{Slutsky} proved that
Frobenius integrability conditions are necessary for the
rationalization of the demand functions (and the inverse demand
functions). On one hand, small perturbations of the demand
functions should violate rationalization conditions. On the other
hand, there was an allusion to the second thermodynamics principle
which has been formulated by Caratheodory also in the form of the
Frobenius integrability conditions. The original question about the
interpretation and satisfiability of the integrability conditions of
the differential form  recovered from the demand functions has been
recognized as a problem that was investigated by such economists as
P.Samuelson, K.Arrow, H.Hauthekker, L.Hurwitz and others. The Frobenius
integrability conditions were reformulated in terms of the strong
axiom of the revealed preference theory, a discrete analog of the
Caratheodory- Rashevsky- Chow criterion. The axiom allows an explicit
testing of the market statistics data (the list of values of the
demand functions at the prescribed points). Computational experiments
showed that violation of the strong axiom took place in the course of
the large-scaled structural changes in the economy similar to the Great
Depession of the thirties. A well-known Afriat- Varian theorem form the
theoretical ground for these experiments. According to one variant of
this theorem the market statistics can be extended to the demand
functions that are rationalized in the class of positive
homogeneous utility functions if and only if it satisfies the
homogeneous strong axiom.

However, the extension figuring in the
Afriat- Varian theorem is not smooth, although smoothness of the
utility functions is an implicit informal assumption. {\em In this
paper we study is it possible to make the extension involved continuous
and additionally to have continuos inverses.} This question in turn is
equivalent to rationalization of the demand functions in the class
of smooth economic indices (positively homogeneous utility functions
and price indices). We note in passing that rationalization of market
statistics in the class of smooth utility functions and
in the class of positive homogeneous utility functions is considered in
\cite{Chiappori}, and  in
\cite{Volosh}, respectively.

\section{Notation}

Denote by $\F$ the set all {\em nonnegative, positive homogeneous,
concave} functions in the nonnegative orthant $\RR_+^n$, i.e. any $f
\in \F$ is a map from $\RR_+^n$ to $\RR_+$ and moreover $f(\cdot)$ is
concave and $f(\lambda x) = \lambda f(x)$ for any nonnegative $\lambda$
and for any $x \in \RR_+^n$.

The Euclidean norm in $\RR^n$ is denoted by $\|\cdot\|$. Let $\BB=\{x
\in \RR^n\,|\, \|x\|\leq1\}$ and $\BB(x_0,r)=\{x\in
\RR^n\,|\,\|x-x_0\|\leq r\}$ be, respectively, the unit Euclidean ball
and the Euclidean ball of the radius $r$ centered at $x_0$.

If $A$ and $B$ are sets in $\RR^n$ their Minkowski sum (i.e. the set of
all points of the form $c=a + b \quad a \in A, b \in B$) is denoted by
$A \oplus B$.

A {\em face} is an intersection of a convex polyhedral set with its
supporting hyperplane. A $(n-1)$-dimensional face of a $n$-dimensional
convex polyhedral set is called a {\em facet}.

The boundary of a set $A \subseteq \RR^n$ is denoted by $Bd(A)$.

Let $A \subseteq \RR^n$ be a closed {\em convex} set and let $x \in
Bd(A)$. Recall that by definition the {\em supporting cone} $T_A(x)$ of
$A$ at $x$ is the intersection of all closed halfspaces containing $A$
and $x$. Respectively, the conjugate (polar) cone $N_A(x)=(T_A(X))^* =
\{p \in \RR^n \,|\, pv \leq 0 \; \forall v \in T_A(x)\}$ is called
the {\em normal cone} of $A$ at $x$.

The boundary of any closed convex body is called a {\em convex
hypersurface} (or a {\em convex curve} when restricted to the plane).

The {\em superdifferential} $\partial f$ of a concave function $f:
\RR^n \rightarrow \RR$ at a point $x_0$ consists of all vectors $p \in
\RR^n$ such that $ p(x-x_0) \geq f(x) - f(x_0)$.

\medskip

\section{Afriat- Varian theorem}

The market statistics
$S = \{p^t, q^t\}, p^t, q^t \in \RR_+^n, t= 1,\ldots
T$, where $p^t$ ¨ $q^t$ are, respectively, the prices and the
volume of purchases in time $t$, is called
{\em rationalizable} in the class of utility functions
$\F$, if there exists such function
$\Q \in \F$ that
$q^t \in \mbox{Argmax }\{\Q(q)\,|\, (p^t q) \leq (p^t q^t), \, q \geq 0
\} \quad t=1,\ldots,T.$

According to the Afriat- Varian theorem (see, e.g.
\cite{Afriat, Shan1}) the following statements are equivalent.
\smallskip

1. The market statistics $S = \{p^t, q^t\}$ is
rationalizable in the class of the utility functions
$\F$.

2. The market statistics $S = \{p^t, q^t\}$ satisfies the
strong homogeneous axiom of the revealed preference
theory, i.e. for any ordered tuple
$\{t(1), \ldots, t(k)\}  \subseteq   \{0, \ldots, T \}$ the inequality
holds:
$$
\begin{array}{l}
(p^{t(1)}q^{t(2)})\, (p^{t(2)}q^{t(3)}) \ldots
(p^{t(k)}q^{t(1)}) \ge \\
\qquad (p^{t(1)}q^{t(1)}) \,(p^{t(2)}q^{t(2)}) \ldots
(p^{t(k)}q^{t(k)}).
\end{array}
$$

3. The following linear system is consistent
\begin{equation}\label{system} \lambda_t (p^t q^t) \leq
\lambda_{\tau} ( p^{\tau} q^t ), \; \lambda_t > 0, \;  t, \tau = 1,
\ldots, T.  \end{equation}
\smallskip

{\em One of the particular utility functions} can be
recovered from an
arbitrary solution of (\ref{system}) by the formula
\begin{equation}\label{utility}
\Q(q) =\min_{i=1, \ldots,t} \lambda_i \, p^i q.
\end{equation}

The superdifferential $\partial \Q(q)$ at any point $q \in
\RR_+^n$
consists of the convex hull of all
(scaled) minimizing prices $\partial \Q(q) = \mbox{conv}_{u \in
U}\{\lambda_u p^u\}$, where $u \in U \subseteq \{1,\ldots, t \}
\Leftrightarrow \lambda_u \, p^u q = \Q(q)$.

Note that the utility function is not {\em uniquely specified} by
the consistency of the system (\ref{system}).

We need the following definition.
\begin{definition}
Let $\Q(\cdot)$ be an arbitrary function from $\F$.
The Lebesque set  $\chi_{\Q}=\{x \in \RR_+^n | P(x) \geq 1\}$ is
called the {\em characteristic set} of $\Q(\cdot)$.
\end{definition}

Denote by $\X$ the set family of all characteristic sets of functions
from the class $\F$.

Next simple lemma summarizes some properties of the characteristic sets.
\begin{lemma}[Characteristic Set]\label{charset}

\begin{enumerate}
\item $\chi \in \X$ iff $\chi$ is a closed convex set in $\RR_+^n$
satisfying hereditary intersection property with any positive ray,
i.e. for any positive ray $r_x=\{\lambda x \,|\, x \in \RR^n, x>0,
\lambda \geq 0\}$ there exists a point $x_0 \in r_x$ such that $r_x
\cap \chi = \{\lambda x_0, \lambda \geq 1\}$.

\item Any $f \in \F$ can be uniquely restored from $\chi_f$ and
vice versa. In other words, there is a canonical bijection
between $\F$ and $\X$.

\item $f \in \F$ is $k$-smooth (i.e. $f \in C^k(\RR_+^n)$) iff $\chi_f$ has
$k$-smooth field of tangential supporting hyperplanes.
\end{enumerate}
\end{lemma}

Let define the {\em gauge transform} of $\Q(\cdot)$:
$$ \P(p)=\inf_{q \geq 0} \frac{qp}{\Q(q)}. $$
As $\Q(\cdot) \in \F$ then the gauge transform $\P(\cdot)$
also belongs to class  $\F$ and forms some dual
index of prices.
By (\cite{Ashm, Fulk})
$$ \Q(q)=\inf_{p \geq 0} \frac{qp}{\P(p)}. $$

Moreover, next lemma shows the correspondence between the level
sets of the dual gauges $\Q(q)$ and $\P(p)$.

\begin{lemma}\label{blocking}
Let  $\Q(q)=1$ (i.e  $q \in Bd(\chi_{\Q}$).
If $p \in \partial \Q(q)$ then $\P(p)=1$ (i.e. $p \in   Bd(\chi_{\P}$ ) and
$q \in \partial \P(p)$ .
(In terminology of \cite{Fulk} the sets $\chi_{\Q}$ and $\chi_{\P}$
form a {\em blocking pair}.) \end{lemma}

{\bf Proof.}  For any $p \in \partial \Q(q)$ and every $\bar q$ the
inequality holds $ p (\bar q - q) \geq \Q(\bar q) - \Q(q).$ Let
$\bar q = \lambda q$, where $\lambda >0$. Then $(\lambda -1) pq \geq
(\lambda -1) \Q(q)$ for every $\lambda >0.$ So we have $pq = \Q(q) =1.$
(The Euler identity for homogeneous but non-smooth $\Q(\cdot).$)
It follows from Kuhn - Tukker
theorem that if $p \in \partial \Q(q)$ then $q \in
Argmax\{\Q(\bar q) | p \bar q \leq p q\}. $ Thus $
\P(p)=\inf_{\bar q \geq 0} \frac{\bar q p}{\Q(\bar q)} = p q = 1.  $ We
have from duality of the gauge transform that $ \Q(q)=\inf_{\bar p \geq
0} \frac{\bar q p}{\P(\bar p)}=1. $ So for every $\bar p \in R^{n}_{+}$
the inequality $q \bar p \geq \P(\bar p)$ holds. Then
for every $\bar p \in R^{n}_{+}$  we obtain $q (\bar p - p) \geq
\P(\bar p) - \P(p),$ i.e. $q \in \partial \P(p)$.

It is shown in \cite{Shan1} that the list of the statements equivalent
to the Afriat- Varian theorem can be enlarged by the following
propositions.

4. Such $\P \in \F$ exists that
$$p^t \in \mbox{Argmax }\{\P(p)\,|\, (q^t p) \leq (q^t p^t), \, p \geq 0 \}
\quad
t=1,\ldots,T.$$
\smallskip

5. The following linear system is consistent
$$
\mu_t (p^t q^t) \leq
\mu_{\tau} ( p^t q^{\tau} ), \, \mu_t > 0,\, t, \tau = 1, \ldots, T.
$$
\smallskip

\begin{definition}
Set $\Q \in \F_k, (k \geq 1)$ if $\Q \in \F$ and additionally both $\Q$
and the dual to
$\Q$ gauge $\P$ belong to $C^k(\RR_+^n)$ (i.e. are $k$-smooth).
\end{definition}
\smallskip

\begin{theorem}
[Smooth Afriat- Varian theorem]

The following statements are equivalent.

\begin{enumerate}

\item The market statistics $S = \{p^t, q^t\}$ is
rationalizable in the class of the utility functions
$\F_k$.

\item\label{twoo}
The market statistics $S = \{p^t, q^t\}$ satisfies the
strong homogeneous axiom of the revealed preference
theory, i.e. for any ordered tuple
$\{t(1), \ldots, t(k)\}  \subseteq   \{0, \ldots, T \}$ the inequality
holds:
$$
\begin{array}{l}
(p^{t(1)}q^{t(2)})\, (p^{t(2)}q^{t(3)}) \ldots
(p^{t(k)}q^{t(1)}) > \\
\qquad (p^{t(1)}q^{t(1)}) \,(p^{t(2)}q^{t(2)}) \ldots (p^{t(k)}q^{t(k)}).
\end{array}
$$

\item\label{three}
The following linear system is (strictly) consistent
\begin{equation}\label{system-strict} \lambda_t (p^t q^t) <
\lambda_{\tau} ( p^{\tau} q^t ), \; \lambda_t > 0, \;  t, \tau = 1,
\ldots, T.  \end{equation}

\item\label{fore}
Such $\P \in \F$ exists that
$$p^t \in \mbox{Argmax }\{\P(p)\,|\, (q^t p) \leq (q^t p^t), \, p \geq 0 \}
\,
t=1,\ldots,T.$$

\item\label{five}
The following linear system is (strictly) consistent
$$
\mu_t (p^t q^t) <
\mu_{\tau} ( p^t q^{\tau} ), \; \mu_t > 0,\; t, \tau = 1, \ldots, T.
$$

\end{enumerate}
\end{theorem}
\smallskip

{\bf Proof.} In fact, it is enough to prove equivalence of the
first and the third statement (all other implications are established
by more or less standard arguments in the framework of the ordinary
Afriat- Varian theorem).

{\bf (1) $\Rightarrow$ (3).}
If market statistics is
rationalizable in the class $\F_k$ then both dual gauges
$\Q(\cdot)$ and $\P(\cdot)$ are smooth and their
supergradients at any points are ordinary gradients (i.e. each
supergradient consists of a single vector only). Therefore, by lemma
\ref{blocking}, the strict solution of the system (\ref{system-strict})
is given by
$\lambda_t = \frac{\|\mbox{grad}(\Q(q^t))\|}{\|p^t\|},\; t=1,\dots,T$.

{\bf (3) $\Rightarrow$ (1).}
Technically, the proof consists in smoothing of the
piecewise- linear utility function that is obtained from some
solution of the system \ref{system-strict}
via suitable convolutions.

Next Lemma is the main technical tool for smoothing.

\begin{definition}
Let $K \subseteq \RR^n$ be a closed convex set and let $n(y), \,
\|n(y)\|=1$ be any unit normal vector at point $y \in Bd(K)$, i.e. the
inclusion $K \subset \{x \in \RR^n \, | \, (n(y),y-x) \leq 0\}$ holds.
Define the {\bf quasi indicator} function of $K$ as follows:
$$\mbox{q-ind}_K(x)=
\left\{
\begin{array}{lr}
1,& x \in K\\
\inf_{y \in Bd(K)} \{1-(n(y),x-y)\},& \mbox{{\em otherwise.}}
\end{array}
\right.
$$
\end{definition}

By construction
$\mbox{q-ind}_K(\cdot)$ is {\em concave}.

\begin{lemma}[Smoothing]\label{smooth}
Let $h(\cdot) \in C^k(\RR^n)$ be any
nonnegative spherically symmetric with respect to the origin
``cap''-function that vanishes outside some $\varepsilon$-ball
$\BB(0,\varepsilon)$ and such that $\int_{\RR^n}{h(x)dx}=1$ (there are
plenty of such ``caps'' even in the $C^{\infty}(\RR^n)$-class).
Let $H= \{x \in \RR^n \, | \, ax \leq 0, \, a \neq 0\}$ be any halfspace
through the origin. Set $\alpha=
\int_{\RR^n}{\mbox{q-ind}_H(x-c) h(c)dc}$.
Then

({\em i})
The convolution
$f = \mbox{q-ind}_K \ast h
= \int_{\RR^n}{\mbox{q-ind}_K(x-c) h(c)dc}$
is a con\-cave func\-tion of the class $C^k(\RR^n)$.

\item\label{middle}

({\em ii})
The following inclusion
holds for the Lebesque set $\chi_{f_{\alpha}}=
\{x\in \RR^n\,|\,f(x)\geq \alpha\} \subseteq K$, in particular,
$\chi_{f_{\alpha}}=K$ if $K$ is any halfspace.

({\em iii})
Let $x \in Bd(K)$ be any boundary point of $K$ and let the boundary
$Bd(K)$ be sufficiently locally flat at $x$, i.e.  $\BB(x,\varepsilon)
\cap K = H_x \cap K$ (here $H_x$ is a supporting halfspace to $K$ at
$x$) then $f(x)=\alpha$.

({\em iv})
If $K \in \X$ is a
characteristic set (of some function from $\F$) then
$\chi_{f_{\alpha}}$ is also a characteristic set from $\X$ (that
corresponds to some function from $\F$).
\end{lemma}
{\bf Proof.}
({\em i})
follows
directly from the properties of convolution (actually it is valid for
an arbitrary locally integrable function).

({\em ii}) and
({\em iii})
are reduced to an easy check.

({\em iv})
follows from
({\em ii})
and from the observation
that for $K \in \X$ all far enough from the origin points of any
positive ray ultimately fall inside $K$ and are sufficiently far
from the boundary of $K$.

It is worth restating here the aforementioned condition for
the existence of the utility function in terms of the
characteristic sets. Namely, the utility function for a given
market statistics $S = \{p^i, q^i\}, p^i, q^i \in \RR_+^n, i= 1 \ldots
t$ exists iff there exists a set $\chi \in \X$, such that
\begin{equation}\label{supergr}
- p^i \in  N_{\chi}(\tilde q^i), \, i=1,\ldots,t.
\end{equation}
Here $\tilde q^i$ is a
(unique) point of intersection of the set $Bd(\chi)$ and the ray
$r_{q^i} = \{\lambda q^i, \lambda \geq 0\}$, $i=1, \ldots, t$ and
$N_{\chi}(\tilde q^i)$ is a corresponding normal cone to $\chi$ at
point $\tilde q^i$.  In particular, the consistent system
(\ref{system}) defines a polyhedral characteristic set $\chi$
(actually, the consistency of (\ref{system}) is a restatement of the
condition (\ref{supergr})).  Thus, in order to smooth a polyhedral
utility function it is sufficient to smooth the corresponding
polyhedral characteristic set preserving the inclusions
(\ref{supergr}).

Now recall that $\P(\cdot)$ is smooth if and only if
the characteristic set $\chi_{\Q}$ has smooth boundary without flat
regions, i.e. the set $\chi_{\Q}$ should have strictly convex boundary.

As above let $\Q(\cdot)$ be the polyhedral utility function
corresponding to the strict solution of (\ref{system-strict}) given by
(\ref{utility}).
Let $\chi_{\Q}$ be the {\em polyhedral} characteristic set of
$\Q(\cdot)$.  Smoothing by convolution is not enough for our
purposes as it preserves some flat regions of $\chi_{\Q}$. To
overcome this difficulty we will slightly change polyhedral
$\chi_{\P}$ replacing its ``flat'' facets by curved ``spherical''
facets preserving the supergradient inclusions (\ref{supergr}).

By construction any facet of $\chi_{\Q}$ is intersected (in its
relative interior) by {\bf a unique} ray
$r_{q^i}=\{x \in \RR_+^n | x = \lambda q^i, \lambda \geq
0\} \; i \in \{1,\dots,t\}$. Let denote the corresponding intersection
points by $s_i, \; i=1,\dots,t$. Now take the system of balls with
equal radii $B_i=\BB(s_i+\rho p_i, \rho), \; i=1,\dots,t$. The radius
$\rho$
should be chosen large enough to satisfy the following
conditions:
\begin{enumerate}
\item any ball $B_i$ is intersected by any nonnegative ray;

\item  for all $i=1,\dots,t$ all points $s_j,\, i \ne j, j=1,\dots,t$
should fall {\em inside} the ball $B_i$.
\end{enumerate}

Set $\hat \chi_{\Q}=
(\cap_{i=1,\dots,t}B_i \cap \chi_{\Q}) \oplus \RR_+^n$.

Namely, to choose $\varepsilon$ we require that for all
$i=1,\dots,t \quad \BB(s_i,\varepsilon) \cap B_i = \BB(s_i,\varepsilon)
\cap \hat \chi_{\Q}$.

Generally, call a boundary point $x$ of $\hat \chi_{\Q}$ belonging to
the boundary of some $B_i$ ${}$ $\varepsilon$-{\em round} if
$\BB(x,\varepsilon) \cap B_i = \BB(x,\varepsilon) \cap \hat \chi_{\Q}$
(thus our requirement for the smallness of $\varepsilon$ means that
all points $s_i$ should be $\varepsilon$-round).

Let $\beta=\mbox{q-ind}_{B_i} \ast
h(s_i+x)$. By our assumptions $\beta$ is a constant for all $i$.
Moreover, by construction this equality holds if we take any point from
the set $U_{s_i}$ of all boundary points of $\hat \chi_{\Q}$ in some
neighborhood of $s_i$ (as all points from $U_{s_i}$ are
$\varepsilon$-round).

As above,
set $\tilde
\chi_{\Q}=\chi_{\varphi_{\beta}}=\{x\in \RR^n\,|\,\varphi (x)\geq
\beta\}$, where $\varphi(\cdot)$ is a convolution of the quasi indicator
function of $\hat \chi_{\Q}$ and $h(\cdot)$. By the reasoning
analogous to the smoothing lemma $\tilde \chi_{\Q} \subset \hat \chi_{\Q}$ and
$\tilde \chi_{\Q} \in \X$. Moreover, for all  $i=1,\dots,t$ all
($\varepsilon$-round) points from $U_{s_i}$ (including $s_i$) belong to
the boundary of $\tilde \chi_{\Q}$ so that supergradient inclusions
(\ref{supergr}) hold for $\tilde \chi_{\Q}$.

Thus
smooth by lemma \ref{charset} index of goods $\tilde \Q(\cdot)$
corresponding to $\tilde \chi_{\Q}$ has a smooth dual index of prices
$\P(\cdot)$ both consistent with the given statistics.

Note that we have proved only $C^1$-smoothness of the indices involved.
$C^k$-smoothness is established by the following argument.

Denote by $\tilde \chi_{\P}$ the characteristic set of the smooth by
construction dual index of prices $\P(\cdot)$. Lemma \ref{blocking}
states that the image of $C^k$-smooth boundary of the characteristic
set $Bd(\tilde \chi_{\Q})$ is mapped under the gradient map $q
\rightarrow \partial \Q$ into the boundary $Bd(\tilde \chi_{\P})$ of the
characteristic set of the corresponding dual (and $C^1$-smooth by
construction) index of prices. (And moreover this map is one-to-one.)
Let prove that the map $q \rightarrow \partial \Q$ has
nonzero Jacobian on the boundary surface $Bd(\tilde \chi_{\Q})$. Then
$C^k$-smoothness of $\P(\cdot)$  immediately follows from the implicit
function theorem. At first, note that
the corresponding Weingarten map is invertible on the
convex surface $Bd(\tilde \chi_{\Q}))$ (Weingarten map surely is
invertible on the spherical patches and hence on the surface obtained
after applying convolution with nonnegative function, i.e. on
$Bd(\tilde \chi_{\Q})$ itself.)  Moreover, it follows from the Euler
identity $q\, \partial \Q=1$ that restricted to the surface $Bd(\tilde
\chi_{\Q})$ our map $q \rightarrow \partial \Q$ can be obtained from
Weingarten map by some smooth scaling of the normal vector.

\begin{rK}
 $C^1$-smoothness of the resulting utility function could be obtained
by a simple geometric construction.
Set $\tilde \chi = (\chi_{\Q} \oplus \varepsilon \BB)
\cap\RR_+^n$, where positive $\varepsilon$ is small enough.
By construction $\tilde \chi \in \X$ and moreover
$\tilde \chi$ has $C^1$-smooth boundary.

The last assertion could be proved as follows.
The set $K \oplus \varepsilon
\BB$ is usually called the {\em outer parallel set} of a convex set $K$
\cite{Leicht}.  It consists of all points within Euclidean distance
$\varepsilon$ from the set $K$. Now the fact that
$\varepsilon$-neighborhood of any closed convex set has $C^1$-smooth
boundary surely belongs to the public mind but we were unable to find
it in the literature. {\bf We'll be grateful for pointing any relevant
references}. For completeness sake we sketch a short geometric proof of
the statement, which was independently communicated to us by M.Arslanov
and S.Chukanov.  Namely, any point $x \in Bd(K) \oplus \varepsilon \BB$
is also a boundary point of a ball $\BB(y,\varepsilon) \subseteq K$
centered at some (boundary) point $y \in K$.  The supporting hyperplane
to $K \oplus \varepsilon \BB$ at $x$ coincides with the supporting
hyperplane to $\BB(y,\varepsilon)$ at $x$ and is thus unique. Hence,
the boundary of $K \oplus \varepsilon \BB$ is $C^1$-smooth.
\end{rK}

\section*{Acknowledgements.} We are grateful to A.Kitaev for proposing
to use convolutions in smoothing procedures and to M.Arslanov and
S.Chukanov for communicating a short proof of the smoothness of an
outer parallel convex body.

\end{document}